\documentclass[12pt,reqno]{amsart}

\usepackage{latexsym}
\usepackage{amssymb}

\renewcommand{\epsilon}{\varepsilon}

\newcommand{\loc}{{\operatorname{loc}}}

\newcommand{\szego}{Szeg\"o }

\newcommand{\kahler}{K\"ahler }

\newcommand{\wt}{\widetilde}

\newcommand{\PP}{{\mathbb P}}
\newcommand{\N}{{\mathbb N}}
\newcommand{\R}{{\mathbb R}}
\newcommand{\C}{{\mathbb C}}

\newcommand{\Z}{{\mathbb Z}}

\newcommand{\CP}{\C\PP}
\renewcommand{\d}{\partial}
\newcommand{\dbar}{\bar\partial}
\newcommand{\ddbar}{\partial\dbar}

\newcommand{\E}{{\mathbf E}}

\newcommand{\half}{{\frac{1}{2}}}

\newcommand{\SU}{{\operatorname{SU}}}

\renewcommand{\phi}{\varphi}

\newcommand{\dcal}{\mathcal{D}}

\newcommand{\fcal}{\mathcal{F}}

\newcommand{\lcal}{\mathcal{L}}

\newcommand{\ocal}{\mathcal{O}}
\newcommand{\pcal}{\mathcal{P}}

\newcommand{\ga}{\gamma}

\newcommand{\la}{\lambda}
\newcommand{\ep}{\varepsilon}

\newcommand{\om}{\omega}
\newcommand{\Om}{\Omega}

\newtheorem{theo}{{\sc Theorem}}[section]

\newtheorem{lem}[theo]{{\sc Lemma}}
\newtheorem{prop}[theo]{{\sc Proposition}}
\newtheorem{example}[theo]{{\sc Example}\rm }
\newenvironment{rem}{\medskip\noindent{\it Remark:\/} }{\medskip}

\title
{Zeros of random polynomials on $\C^m$}

\author{Thomas Bloom}
\address{Department of Mathematics,
University of Toronto,
Toronto, ON,
CANADA M5S 3G3}
\email {bloom@math.toronto.edu}

\author{Bernard Shiffman}
\address{Department of Mathematics, Johns Hopkins University, Baltimore, MD
21218, USA} \email{shiffman@math.jhu.edu}

\thanks{Research partially supported by NSERC (first author) and NSF grant
DMS-0100474 (second author).}

\date{May 27, 2006}

\begin{document}

\begin{abstract} For a regular compact set $K$ in
$\mathbb{C}^m$ and a measure $\mu$ on $K$ satisfying the
Bernstein-Markov inequality, we consider the ensemble
$\mathcal{P}_N$ of polynomials of degree $N$, endowed with the
Gaussian probability measure induced by
$L^2(\mu)$.  We
show that for large $N$, the simultaneous
zeros of $m$ polynomials in $\mathcal{P}_N$ tend to concentrate
around the Silov boundary of $K$; more precisely, their expected
distribution  is asymptotic to $N^m \mu_{eq}$, where
$\mu_{eq}$ is the equilibrium measure of $K$. For the
case where $K$ is the unit ball, we give scaling
asymptotics for the expected distribution of zeros as
$N\to\infty$.\end{abstract}

\maketitle

 \section{Introduction}

A classical result due to Hammersley \cite{Ha} (see also \cite{SV}),
loosely stated, is that the zeros of a random
complex polynomial
\begin{equation}
f(z)=\sum^N_{j=0} c_j z^j
\end{equation}
mostly tend towards the unit circle $|z|=1$ as the degree
$N\to\infty$, when the coefficients $c_j$ are independent complex
Gaussian random variables of mean zero and variance one.
In this paper, we will prove a
multivariable result (Theorem \ref{general}), a special
case (Example \ref{ex-poly}) of which shows, loosely stated,
that the common zeros of $m$ random complex
polynomials in $\C^m$,
\begin{equation}
f_k (z)=\sum_{|J|\le N} c_J^k\, z_1^{j_1}\cdots z_m^{j_m}\qquad
{\rm for}\ k=1,\dots, m\;,
\end{equation}
tend to concentrate on the product of the unit circles $|z_j|=1$
($j=1,\dots, m$) as $N\to\infty$,  when the coefficients $c^k_J$
are i.i.d.\ complex Gaussian random variables.

The following is our basic setting:  We let $K$ be a compact set
in $\C^m$ and let $\mu$ be a Borel probability measure on $K$.  We
assume that $K$ is non-pluripolar and we let $V_K$ be its
pluricomplex Green function.  We also assume that $K$ is regular
(i.e., $V_K=V_K^*$) and that $\mu$ satisfies the Bernstein-Markov
inequality (see \S\ref{background}).  We give the space $\pcal_N$
of holomorphic polynomials of degree $\le N$ on $\C^m$ the
Gaussian probability measure $\ga_N$ induced by the Hermitian
inner product
\begin{equation}\label{inner}(f,g) = \int_Kf \bar g\,d\mu\;.\end{equation}
The Gaussian measure $\ga_N$ can be described as follows:  We write
$f=\sum_{j=1}^{d(N)}c_jp_j$, where $\{p_j\}$ is  an orthonormal
basis of $\pcal_N$ with respect to \eqref{inner} and $d(N)=\dim
\pcal_N={{N+m}\choose m}$.  Identifying $f\in\pcal_N$ with
$c=(c_1,\dots,c_{d(N)})\in \C^{d(N)}$, we have
\begin{equation*}d\ga_N(s)=\frac{1}{\pi^{d(N)}}e^
{-|c|^2}dc\,.\end{equation*}
(The measure $\ga_N$  is independent of the choice of orthonormal basis $\{p_j\}$.)
In other words, a random polynomial in
the ensemble $(\pcal_N,\ga_N)$ is a polynomial $f = \sum_j c_j\,p_j$, where the
$c_j$ are independent complex Gaussian random variables with mean 0 and variance 1.

Our main result, Theorem \ref{general}, gives  asymptotics for the
expected zero current of $k$ i.i.d.\ random polynomials ($1\le
k\le m$).  In particular,  the expected
distribution $\E( Z_{f_1,\dots,f_m})$ of simultaneous zeros of
$m$ independent random polynomials in $(\pcal_N,\ga_N)$ has the asymptotics
\begin{equation}\label{equilib}\frac 1 {N^m} \, \E( Z_{f_1,\dots,f_m}) \to
\mu_{eq}\quad  weak^*\;,\end{equation} where $\mu_{eq}= (\frac
i\pi \ddbar V_K)^m$ is the equilibrium measure of $K$.  Here,
$\E(X)$ denotes the expected value of a random variable $X$.

The reader may notice from \eqref{equilib} that the distributions
of zeros for the measures on $\pcal_N$ considered here are quite
different from those of the $\SU(m+1)$ ensembles studied, for
example, in \cite{SZ,SZvar,BSZ1,BSZ2,DS}. The Gaussian measure on
the $\SU(m+1)$ polynomials is based on the inner product $$\langle
f,g\rangle_N = \int_{S^{2m+1}}F_N \;\overline{G_N},$$ where $
F_N,G_N \in \C[z_0,z_1,\dots,z_m]$ denote the degree $N$
homogenizations of $f$ and $g$ respectively. It follows easily
from the $\SU(m+1)$-invariance of the inner product that the
expected distribution of simultaneous zeros equals $\frac
{N^m}{\pi^m} \om^m$ (exactly), where $\om$ is the Fubini-Study
\kahler form (on $\C^m\subset \CP^m$).  We note that, unlike
\eqref{inner}, this inner product depends on $N$; indeed,
$\|z^J\|_N^2 = \frac{m!(N-|J|)!j_1!\cdots j_m!}{(N+m)!}$
\cite[(30)]{SZ}.

%%% INSERT B %%%
In this paper, we  also give scaling limits for the expected zero
density in the case of the unit ball in $\C^m$  (Theorem
\ref{sphere-asymp}). The problem of finding scaling limits for
more general sets in $\C^m$ remains open. Another open problem is
to establish the multivariable version of the following one
variable result: For a regular subset $K\subset \C$,  it is known
(see \cite{SZ,Bl2}) that with probability one, a sequence
$\{f_N\}_{N=1, 2,\dots}$ of random polynomials of increasing
degree satisfies:
\[
\lim\limits_{N\to \infty} {1\over N} Z_{f_N}=\mu_{eq} \quad
weak^*\;.
\]

\section{Background}\label{background}

%%%% INSERT C  %%%%
We let $\lcal$ denote the Lelong class of
plurisubharmonic (PSH) functions on $\C^m$ of
at most logarithmic growth at $\infty$. That is
\begin{equation}
\lcal:=\{u\in \mbox{PSH} (\C^m)\mid u(z)\le \log^+
\|z\| +O(1)\}
\end{equation}
For $K$ a compact subset of $\C^m$, we define its pluricomplex
Green function $V_K(z)$ via
\begin{equation}
V_K(z)=\sup\{u (z)\mid u \in \lcal,\ u\le 0\ \ {\rm on}\ \ K\}.
\end{equation}
We will assume $K$ is regular, that is by
definition, $V_K$ is continuous on $\C^m$ (and so
$V_K=V_K^*$, its uppersemicontinuous regularization). The function $V_K$ is
a locally bounded PSH function on $\C^m$ and, in fact
 \begin{equation}\label{Vk}
 V_K - \log^+ \|z\|=O (1)\;.
 \end{equation}
 By a basic result of Bedford and Taylor \cite{BT1} (see \cite{Klimek}), the
 complex Monge-Amp\`ere operator $(dd^c)^m =(2i\partial\bar\partial)^m$
 is defined on any locally bounded PSH function $\C^m$ and
 in particular on $V_K$. The equilibrium measure of $K$
 is defined by (see \cite[Cor.~5.5.3]{Klimek})
 \begin{equation}
 \mu_{eq} (K):=\left({i\over \pi} \partial\bar\partial V_K\right)^m
 \end{equation}

 Since $V_K$ satisfies \eqref{Vk}, it is a positive Borel measure,
 here normalized to have mass 1. The support of the
 measure $\mu_{eq} (K)$ is the Silov boundary of $K$ for
 the algebra of entire analytic functions \cite{BT2}.
In one variable, i.e.\ $K\subset \C$, $V_K$ is the Green
 function of the unbounded component of $\C \setminus K$ with
 a logarithmic pole at $\infty$, and $\mu_{eq} (K)={1\over 2\pi}
 \Delta V_K$, where $\Delta$ is the Laplacian \cite{Ra}.

 Let $\mu$ be a finite positive Borel measure on $K$. The measure
 $\mu$ is said to satisfy a Bernstein-Markov (BM)
 inequality, if, for each $\varepsilon >0$ there is a
 constant $C=C(\varepsilon)>0$ such that
\begin{equation}\label{BM}
\|p\|_K\le C e^{\ep\deg(p)} \|p\|_{L^2(\mu)}
\end{equation}
for all holomorphic polynomials $p$. Essentially, the BM
inequality says that the $L^2$ norms and the sup norms of a
sequence of holomorphic polynomials of increasing degrees are
``asymptotically equivalent".

The question arises as to which measures actually satisfy the BM
inequality. It is a result of Nguyen-Zeriahi \cite{NZ} combined
with \cite[Cor.~5.6.7]{Klimek} that for $K$ regular, $\mu_{eq}(K)$
satisfies BM. This fact is used in Examples
\ref{ex-poly}--\ref{ex-ball}. In \cite[Theorem 2.2]{Bl1}, a
``mass-density" condition for a measure to satisfy BM was given.
(See also \cite{BL}.)

Our proof uses the {\it probabilistic Poincar\'e-Lelong formula\/}
for the zeros of random functions (Proposition \ref{PL} below).
Considering a slightly more general situation, we let $g_1,\dots,
g_d$ be holomorphic functions with no common zeros on a domain
$U\subset \C^m$.  (We are interested in the case where $U=\C^m$
and  $\{g_j\}$ is an orthonormal basis of $\pcal_N$ with respect
to the inner product \eqref{inner}, as discussed above.) We let
$\fcal$ denote the ensemble of random holomorphic functions of the
form $f = \sum c_j\,g_j$, where the $c_j$ are independent complex
Gaussian random variables with mean 0 and variance 1. We consider
the {\it \szego kernel}
$$S_\fcal(z,w) = \sum_{j=1}^d g_j(z)\,\overline{g_j(w)}\;.$$ For
the case where the $g_j$ are orthonormal functions with respect to
an inner product on $\ocal(U)$, $S_\fcal(z,w)$ is the kernel for
the orthogonal projection onto the span of the $g_j$.

Under the assumption that the $g_j$ have no common zeros, it is
easily shown using Sard's theorem (or a variation of Bertini's
theorem) that for almost all $f_1,\dots,f_k\in\fcal$, the
differentials $df_1,\dots,df_k$ are linearly independent at all
points of the zero set
$$\loc(f_1,\dots,f_k):=\{z\in U:f_1(z)=\cdots=f_k(z)=0\}\;.$$
This condition implies that the complex hypersurfaces $\loc(f_j)$
are smooth and intersect transversely, and hence
$\loc(f_1,\dots,f_k)$ is a codimension $k$ complex submanifold of
$U$.  We then let $Z_{f_1,\dots,f_k}\in\dcal'{}^{k,k}(U)$ denote
the current of integration over $\loc(f_1,\dots,f_k)$:
$$\big(Z_{f_1,\dots,f_k},
\phi\big) = \int_{\loc(f_1,\dots,f_k)}\phi\;, \qquad \phi\in
\dcal^{m-k,m-k}(U)\;.$$
 We shall use the following Poincar\'e-Lelong
formula from \cite{SZnewton, SZvar}:

\begin{prop}\label{PL} The expected zero current of $k$ independent random functions
$f_1,\dots,f_k\in \fcal$ is given by
$$\E(Z_{f_1,\dots,f_k})= \left(\frac i {2\pi} \ddbar \log S_\fcal(z,z)\right)^k.$$
\end{prop}
The proof follows by a verbatim repetition of the proof of
Proposition 5.1 in \cite{SZnewton} (which gives the case where the
$g_j$ are normalized monomials with exponents in a Newton
polytope). The codimension $k=1$ case was given in \cite{SZ} (for
sections of holomorphic line bundles), and in dimension 1 by
Edelman-Kostlan \cite{EK}. (The formula also holds for
infinite-dimensional ensembles; see \cite{So, SZvar}.)  The
general case follows from the codimension 1 case together with the
fact that
\begin{equation}\label{wedge}\E(Z_{f_1,\dots,f_k})=\E(Z_{f_1})\wedge\cdots
\wedge \E(Z_{f_k})= \E(Z_{f})^k\;,\end{equation} which is a
consequence of the independence of the $f_j$. The wedge product of
currents is not always defined, but $Z_{f_1}\wedge\cdots \wedge
Z_{f_k}$ is almost always defined (and equals $Z_{f_1,\dots,f_k}$
whenever the hypersurfaces $\loc(f_j)$ are smooth and intersect
transversely), and a short argument given in \cite{SZnewton}
yields \eqref{wedge}. (In fact, the left equality of \eqref{wedge}
holds for independent non-identically-distributed $f_j$, as proven
in \cite{SZnewton}.)  We note that the expectations in
\eqref{wedge} are smooth forms.

\section{Random polynomials on polynomially convex sets}

\begin{theo}\label{general} Let $\mu$ be a Borel probability
measure on a regular compact set $K\subset
\C^m$, and suppose that $(K,\mu)$ satisfies the Bernstein-Markov
inequality. Let $1\le k\le m$, and let $(\pcal_N^k, \ga_N^k) $
denote the ensemble of $k$-tuples of i.i.d.\ Gaussian random
polynomials of degree $\le N$ with the Gaussian measure $d\ga_N$
induced by $L^2(\mu)$. Then
$$\frac 1{N^k} \E_{\ga_N^k}(Z_{f_1,\dots,f_k}) \to \left(\frac i{\pi}\ddbar
V_K\right)^k \qquad weak^*, \quad \mbox{as\ } N\to \infty\;,$$
where $V_K$ is the pluricomplex Green function of $K$ with pole at
infinity.
\end{theo}

To prove Theorem \ref{general}, we consider the \szego kernels
$$S_N(z,w):=S_{(\pcal_N,\ga_N)}(z,w)= \sum_{j=1}^{d(N)}
p_j(z)\overline{p_j(w)}\;,$$ where $\{p_j\}$ is an
$L^2(\mu)$-orthonormal basis for $\pcal_N$.
 Our proof  is based on approximating
the extremal function $V_K$  by the (normalized) logarithms of the
\szego kernels $S_N(z,z)$ (Lemma \ref{szego}).

We begin by considering the polynomial suprema
\begin{equation}\label{PhiN} \Phi_N^K(z)= \sup\{|f(z)|:
{f\in\pcal_N},\ \|f\|_K\le 1\}\;.\end{equation}  Since $\frac
1N\log f\in\lcal$, for $f\in\pcal_N$, it is clear that $\frac 1{N}
\log \Phi_N^K \le V_K$, for all $N$. Pioneering work of Zaharjuta
\cite{Za} and Siciak \cite{Si,Si2} established the convergence of
$\frac 1N \log\Phi_N^K$ to $V_K$. The uniform convergence when $K$
is regular seems not to have been explicitly stated and we give
the proof below.

\begin{lem}\label{uniform} Let $K$ be a regular compact set in
$\C^m$. Then
$$\frac 1{N} \log \Phi_N^K(z) \to V_K(z)$$ uniformly on compact
subsets of $\C^m$. \end{lem}

\begin{proof} We first note that
$1\le \Phi_j\le \Phi_j\Phi_k \le \Phi_{j+k}$, for $j,k \ge0$. By a
result of Siciak \cite{Si} and Zaharjuta \cite{Za} (see
\cite[Theorem 5.1.7]{Klimek}),
\begin{equation}\label{pointwise} V_K(z) = \lim_{N\to\infty}
\frac 1N\log\Phi_N^K(z) = \sup_N \frac
1N\log\Phi_N^K(z)\;,\end{equation} for all $z\in\C^m$.

 We use  the regularity of
$K$ to show that the convergence is uniform: let $$\psi_N = \frac
1N \log \Phi_N^K \ge 0\;.$$  Thus for $N,k\ge 1$, $j\ge 0$, we
have
$$Nk\,\psi_{Nk} +j\,\psi_j \le(Nk+j)\,\psi_{Nk+j}\,.$$ Since
$\psi_N \le \psi_{Nk}$, we then obtain the inequality
\begin{equation}\label{subadd}  \psi_{Nk+j} \ge \frac
{Nk}{Nk+j} \psi_N+ \frac {j}{Nk+j} \psi_j \ge \frac{Nk}{Nk+j}
\psi_N\;.\end{equation}

Fix  $\ep>0$.  For each $a\in\C^m$, we choose $N_a\in\Z^+$ such
that
$$V_K(a)-\psi_{N_a}(a) <\epsilon \qquad \mbox{and}\qquad
\frac{V_K(a)}{N_a} <\ep\;,$$ and then choose a neighborhood $U_a$
of $a$ such that
$$|V_K(z)-V_K(a)|<\ep,\quad
\psi_{N_a}(z)\ge\psi_{N_a}(a)-\ep,\quad \frac{V_K(z)}{N_a}
<\ep,\qquad \mbox{for }\ z\in U_a\;.$$ Now let $N\ge N_a^2$, and
write $N=N_ak+j$, where $k\ge N_a,\ 0\le j<N_a$. By
\eqref{pointwise}--\eqref{subadd}, we have
\begin{equation}\label{sub2}0\le  V_K-\psi_N\le V_K-\frac
{N_ak}{N_ak+j} \psi_{N_a} \le V_K-\frac {N_a}{N_a+1}\psi_{N_a}\le
V_K-\psi_{N_a} +\frac {1}{N_a+1}V_K\;.\end{equation} Hence, for
all $N\ge N_a^2$ and for all $z\in U_a$, we have
\begin{eqnarray}0\ \le \ V_K(z)-\psi_N(z) &<& V_K(z)-\psi_{N_a}(z)
+\ep\nonumber\\ &=&[V_K(a)-\psi_{N_a}(a)] +[V_K(z)-V_K(a)]+
[\psi_{N_a}(a)-\psi_{N_a}(z)]+\ep\nonumber \\&< &
4\ep\,.\label{sub3}
\end{eqnarray}
Hence for each compact  $A\subset\C^m$, we can cover $A$ with
finitely many $U_{a_i}$, so that we have by
\eqref{sub3},
$$\|V_K-\psi_N\|_A \le 4\ep \qquad \forall\ N\ge \max_i N_{a_i}^2\;.$$
\end{proof}

\begin{lem}\label{szego-est}
For all $\ep >0$, there exists $C=C_\ep>0$ such that
$$\frac {1}{d(N)} \le \frac {S_N(z,z)}{\Phi_N^K(z)^2}\le
C\,e^{\ep N} d(N).$$
\end{lem}

\begin{proof}  Let $f\in\pcal_N$ with $\|f\|_K\le 1$. Then
\begin{eqnarray*} |f(z)| &=& \left|\int_K S_N(z,w)f(w)\,d\mu(w)\right| \ \le\
\int_K|S_N(z,w)|\,d\mu(w)\\ &\le & \int_KS_N(z,z)^\half S_N(w,w)^\half \,d\mu(w)\ =\
S_N(z,z)^\half\, \|S_N(w,w)^\half\|_{L^1(\mu)}\\
&\le & S_N(z,z)^\half\, \|1\|_{L^2(\mu)} \,
\|S_N(w,w)^\half\|_{L^2(\mu)}\ = \  S_N(z,z)^\half\, d(N)^\half
\;.\end{eqnarray*} Taking the supremum over $f\in\pcal_N$ with
$\|f\|_K\le 1$, we obtain the left inequality of the lemma.

To verify the right inequality, we let $\{p_j\}$  be a sequence of
$L^2(\mu)$-orthonormal polynomials, obtained by applying
Gram-Schmid to a sequence of monomials of non-decreasing degree,
so that $\{p_1,\dots p_{d(N)}\}$ is an orthonormal  basis of
$\pcal_N$ (for each $N\in\Z^+$). By the Bernstein-Markov
inequality
\eqref{BM}, we have
$$\|p_j\|_K \le C\,e^{\ep\,\deg p_j}$$ and hence $$|p_j(z)|\le
\|p_j\|_K\,\Phi^K_{\deg p_j}(z) \le  C\,e^{\ep\,\deg
p_j}\,\Phi^K_{\deg p_j}(z) \le C\,e^{\ep N}\,\Phi_N^K(z)\;,
\quad\mbox{for }\ j\le d(N).$$  Therefore,
$$S_N(z,z) = \sum_{j=1}^{d(N)} |p_j(z)|^2 \le
d(N)\,C^2\,e^{2\ep N} \Phi_N^K(z)^2\;.$$
\end{proof}

\begin{lem}\label{szego} Under the hypotheses of Theorem \ref{general}, we have
$$\frac 1{2N} \log S_N(z,z) \to V_K(z)$$ uniformly on compact
subsets of $\C^m$. \end{lem}

\begin{proof} Let $\ep>0$ be arbitrary. Recalling that
$d(N)={N+m \choose m}$, we have by Lemma \ref{szego-est},
$$-\frac mN \log(N+m) \le \frac 1{N}\log\left(\frac {S_N(z,z)}{\Phi_N^K(z)^2}
\right) \le\frac {\log C}N + \ep + \frac mN \log(N+m)\;.$$ Since
$\ep>0$ is arbitrary, we then have
\begin{equation}\label{SNPhiN}\frac 1{N}\log\left(\frac
{S_N(z,z)}{\Phi_N^K(z)^2} \right) \to 0\;.\end{equation} The
conclusion follows from Lemma \ref{uniform} and \eqref{SNPhiN}.
\end{proof}

\medskip\noindent {\it Proof of Theorem \ref{general}:\/}  It follows from Lemma
\ref{szego} and the fact that the complex Monge-\`Ampere operator
is continuous under uniform limits \cite{BT1},
 $$\left(\frac i{2\pi N} \ddbar\log S_N(z,z)\right)^k \to
\left( \frac i{\pi} \ddbar V_K(z)\right)^k\qquad weak^*\,.$$ The
conclusion then follows from Proposition \ref{PL}.\qed

%%%% INSERT D  %%%%%
\begin{example}\label{ex-poly}\rm Let $K$ be the unit  polydisk in $\C^m$. Then
$V_K= \max^m_{j=1} \log^+|z_j|$, the Silov boundary of $K$ is the
product of the circles $|z_j|=1$ ($j=1, \dots, m$), and
$d\mu_{eq}=({1\over 2\pi})^m d\theta_1\cdots d\theta_m$ where
$d\theta_j$ is the angular measure on the circle $|z_j|=1$.

The monomials $z^J:=z_1^{j_1} \cdots z^{jm}_m$, for $|J|\le N$,
form an orthonormal basis for $\pcal_N$. A random polynomial in
the ensemble is of the form
\[
f(z)=\sum_{|J| \le N} c_J z^J
\]
where the $c_J$ are independent complex Gaussian random variables
of mean zero and variance one. By Theorem \ref{general},
$\E_{\ga_N^k}(Z_{f_1 ,\dots, f_m}) \to ({1\over 2\pi})^m
d\theta_1\cdots d\theta_m$ weak$^*$,  as $N\to \infty$. In
particular, the common zeros of $m$ random polynomials tend to the
product of the unit circles $|z_j|=1$ for $j=1,\dots,
m$.\end{example}

\begin{example}\label{ex-ball}\rm Let $K$ be the unit ball $\{
\|z\|\le 1\}$ in $\C^m$. Then the Silov boundary of $K$ is its
topological boundary $\{ \|z\|=1\}$, $V_K(z)=\log^+ \|z\|$, and
$\mu_{eq}$ is the invariant hypersurface measure on $\|z\|=1$
normalized to have total mass one. \end{example}
%%%%%%%%%%%%

\section{Scaling limit zero density for orthogonal polynomials on
$S^{2m-1}$}\label{sphere}

%%% INSERT D  %%%%
Examples \ref{ex-poly} and \ref{ex-ball} both reduce to the unit
disk in the one variable case. In that case, detailed scaling
limits are known (see, for example, \cite{IZ}). For a more general
compact set $K\subset \C$ with an analytic boundary, scaling
limits are found in \cite{SZequil}.
%%%%

In this section, we consider the case where $K=\{z\in \C^m:
\|z\|\le 1\}$ is the unit ball and $\mu$ is its equilibrium
measure, i.e.\ invariant measure on the unit sphere $S^{2m-1}$. We
have the following scaling asymptotics for the expected
distribution of zeros of $m$ random polynomials orthonormalized on
the sphere:

\begin{theo}\label{sphere-asymp} Let $(\pcal_N^m, \ga_N^m) $
denote the ensemble of $m$-tuples of i.i.d.\ Gaussian random
polynomials of degree $\le N$ with the Gaussian measure $d\ga_N$
induced by $L^2(S^{2m-1},\mu)$, where $\mu$ is the invariant
measure on the unit sphere $S^{2m-1}\subset \C^m$. Then
\begin{equation*}
\E_{\ga_N^m}(Z_{f_1,\dots,f_m}) = D_N\left(\log{\|z\|^2}\right)\,
\left(\frac i{2} \ddbar \|z\|^2\right)^m,\end{equation*} where
$$\frac
1{N^{m+1}} D_N\left(\frac uN\right) =  \frac {1}{\pi^m}\,  F''_m(u)\,F'_m(u)^{m-1}
+O\left(\frac 1N\right)\;,\quad F_m(u) = \log \left[ \frac
{d^{m-1}}{du^{m-1}}\left(\frac {e^u-1}u\right)\right].
$$
\end{theo}

\begin{proof} We write $$z^J=z_1^{j_1}\cdots z_m^{j_m}\,,\qquad z=(z_1,\dots,z_m),\
J=(j_1,\dots, j_m)\;.$$  An easy computation yields
\begin{equation}\label{ZJ} \int_{S^{2m-1}}|z^J|^2\,d\mu(z) = \frac{(m-1)!j_1!\cdots
j_m!}{(|J|+m-1)!}  = \frac 1 {{|J|+m-1\choose m-1}{|J|\choose
J}}\;,\end{equation} where $$|J|=j_1+\cdots +j_m\,,\qquad
{|J|\choose J}= \frac {|J|!}{j_1!\cdots j_m!}\;.$$  Thus an
orthonormal basis for $\pcal_N$ on $S^{2m-1}$ is:
\begin{equation}\label{ortho} \phi_J(z) = {|J|+m-1\choose m-1}^\half {|J|\choose
J}^\half \,z^J\;, \qquad |J|\le N\;.\end{equation}

We have
\begin{eqnarray*} S_N(z,z) = \sum_{|J|\le N} |\phi_J(z)|^2 &=&
\sum_{k=0}^N{k+m-1\choose m-1} \sum_{|J|=k} {k\choose J}|z_1|^{2j_1} \cdots
|z_m|^{2j_m}
\\&=&\sum_{k=0}^N{k+m-1\choose m-1}
\|z\|^{2k}\,.\end{eqnarray*}  Hence
\begin{equation}
\label{SNgN}S_N(z,z)=g_N(\|z\|^2)\;,\quad\mbox{where } \ g_N(x)=
\sum_{k=0}^N{k+m-1\choose m-1}\; x^k\;.\end{equation}  We note
that
$$g_N = \frac 1{(m-1)!}\,G_N^{(m-1)} \;,\quad\mbox{where } \
G_N(x)=\frac{1-x^{N+m}}{1-x}\;.$$

We denote by $O(\frac 1N)$ any function
$\la(N,u)=\la_N(u):\Z^+\times \R\to\R$ satisfying:
\begin{equation}\label{oN} \forall R>0,\ \forall j\in \N, \
\exists C_{Rj}\in\R^+ \mbox{ \ \ such that } \sup_{|u|<R}|\la_N^{(j)}(u)| <
\frac {C_{Rj}}{N}\;.\end{equation}  We note that
$$N\log \left(1+\frac uN\right) = u + u^2\,O\left(\frac 1N\right) \qquad(\mbox{for }\
|u|<N)\;,$$ and hence $$\left(1+\frac uN\right)^N= e^u + u^2 \,O\left(\frac
1N\right)\;.$$  Thus we have

\begin{equation}\label{scaling} \frac 1N \,G_N\left(1+\frac uN\right) = \frac {e^u-1}u
+ O\left(\frac 1N\right)\;.
\end{equation}  Hence
\begin{equation}\label{scaling1} \frac 1{N^m} \;g_N\left(1+\frac uN\right) = \frac
1{(m-1)!}\;\frac {d^{m-1}}{du^{m-1}}\left(\frac {e^u-1}u\right) + O\left(\frac
1N\right)\;.
\end{equation}
Therefore
\begin{equation}\label{FNinfinity}\log \left[ \frac
{(m-1)!}{N^m} \,g_N\left(1+\frac uN\right) \right]= F_m(u) + O\left(\frac
1N\right) \;,\end{equation} where $F_m$ is given in the statement
of the theorem.

Since the zero distribution is invariant under the SO$(2m)$-action
on $\C^m$, we can write \begin{equation}\label{def}
\E_{\ga_N^m}(Z_{f_1,\dots,f_m}) = D_N\left(\log{\|z\|^2}\right)\,
\left(\frac i{2} \ddbar \|z\|^2\right)^m.\end{equation} Then
$D_N(\frac uN)$ is the density at the point
$$z^N:=\left(\frac 1{\sqrt m} \,e^{u/2N},\dots, \frac 1{\sqrt m}
\,e^{u/2N}\right)\in\C^m\;, \qquad  \|z^N\|^2=e^{u/N}\;.$$

We shall compute using the local coordinates $\zeta_j =
\rho_j+i\theta_j=\log z_j$. Let
$$\Om= \left(\frac i2 \ddbar \sum|\zeta_j|^2\right)^m\;.$$  By Proposition
\ref{PL} and \eqref{SNgN}, we have
\begin{eqnarray} \E_{\ga_N^m}(Z_{f_1,\dots,f_m})=
\left(\frac 1{2\pi}\right)^m \det \left(\frac 12\,
\frac{\d^2}{\d\rho_j\d\rho_k} \log g_N\left(\sum
e^{2\rho_j}\right)\right) \Om\label{EZrho} \;.\end{eqnarray} We
note that \begin{equation}\label{Omega}\Om=m^m\left[1 +
O\left(\frac 1N\right)\right]\left(\frac i2 \ddbar
\|z\|^2\right)^m \quad \mbox{ at the point }\ z^N\;.\end{equation}
We let $\mathbf{1}$ denote the $m\times m$ matrix all of whose
entries are equal to 1 (and we let $I$ denote the $m\times m$
identity matrix). By \eqref{FNinfinity} and
\eqref{EZrho}--\eqref{Omega}, we have
\begin{eqnarray*} D_N\left(\frac uN \right)&=& \left(\frac
{m}{2\pi}\right)^m \left[1 + O\left(\frac 1N\right)\right]\\&&
\quad \times \det \Big( 2\,m^{-2}\,e^{2u/N} (\log
g_N)''(e^{u/N})\,\mathbf{1} +
2\, m^{-1}\,e^{u/N} (\log g_N)'(e^{u/N})\,I\Big)\\
&=& \frac {1}{\pi^m} \left[1 + O\left(\frac
1N\right)\right]\det \Big(m^{-1} N^2\, F_m''(u)
\,\mathbf{1} +N\,F_m'(u)\,I\Big) \;.\end{eqnarray*} Therefore,
\begin{eqnarray*} \frac 1{N^{m+1}}\,D_N\left(\frac uN \right)&=&\frac
1{N^{m+1}\,\pi^m}
 \left[1 + O\left(\frac
1N\right)\right]\\&& \quad \times\left\{
\big[N\,F_m'(u)\big]^{m}+\ m \left[ m^{-1} N^2\,
F_m''(u)\right]\big[N\,F_m'(u)\big]^{m-1}\right\}\\
&=&  \frac {1}{\pi^m}\,F_m''(u)\,F_m'(u)^{m-1}
+O\left(\frac 1N\right)\;.\end{eqnarray*}
\end{proof}

\begin{rem} There is a similarity between the scaling asymptotics
of Theorem \ref{sphere-asymp} and that of the one-dimensional
$\SU(1,1)$ ensembles in \cite{BR} with the norms $\|z^j\|={L-1+j
\choose j}^{-1/2}$, for $L\in\Z^+$.   Then the expected
distribution of zeros of random $\SU(1,1)$ polynomials of degree
$N$ has the asymptotics \cite[Th.~2.1]{BR}:
\begin{equation*}
\E_N(Z_{f}) = \wt D_N\left(\log{|z|^2}\right)\, \frac i{2}
dz\wedge d\bar z\;,\end{equation*} where (in our notation)
$$\frac
1{N^2} \wt D_N\left(\frac uN\right) = \frac {1}{\pi}\,
F''_{L-1}(u) +O\left(\frac 1N\right)\;.
$$

\end{rem}

\end{document}